\documentclass{article}    
\usepackage{amstext}
\def\qed{$\rlap{$\sqcap$}\sqcup$}
\usepackage{latexsym}

\begin{document}           

{\ }\\
\begin{center}
{\huge {\bf When is there a unique socle-vector associated to a given $h$-vector?}} \\ [.250in]
{\large FABRIZIO ZANELLO\\
Dipartimento di Matematica, Universit\`a di Genova, Genova, Italy.\\E-mail: zanello@dima.unige.it}
\end{center}
{\large

{\ }\\
\\
ABSTRACT. First, we construct a bijection between the set of $h$-vectors and the set of socle-vectors of artinian algebras. As a corollary, we find the minimum codimension that an artinian algebra with a given socle-vector can have. Then, we study the main problem in the paper: determining when there is a unique socle-vector for a given $h$-vector. We solve the problem completely if the codimension is at most 3.\\
\\
\section{Introduction}
\indent

The purpose of this paper is to begin the study of the possible socle-vectors of algebras having a given $h$-vector $h$. Our algebras are standard, graded and artinian, and will be denoted by $A=R/I$, where $R=k[x_1,...,x_r]$, $I$ is a homogeneous ideal of $R$ and $k$ is a field of characteristic zero.\\\indent
This article is motivated by the fact that much work has been done on the converse problem, that of determining the $h$-vectors associated to a given socle-vector $s$ (often when $s$ is special, e.g. level or Gorenstein), while the problem we are going to investigate here has never been specifically addressed. Moreover, determining all the socle-vectors of the algebras sharing the same $h$-vector can be seen as part of the more general (and very difficult) problem of determining all the minimal free resolutions (MFR's) associated to a given $h$-vector.\\
\\\indent
We first recall some basic definitions. The {\it $h$-vector} of $A$ is $h(A)=h=(h_0,h_1,...,h_e)$, where $h_i=\dim_k A_i$ and $e$ is the last index such that $\dim_k A_e>0$. Since we may suppose, without loss of generality, that $I$ does not contain non-zero forms of degree 1, $r=h_1$ is defined to be the {\it codimension} of $A$ (or of $h$).\\\indent 
The {\it socle} of $A$ is the annihilator of the maximal homogeneous ideal $\overline{m}=(\overline{x_1},...,\overline{x_r})\subseteq A$, namely $soc(A)=\lbrace a\in A {\ } \mid {\ } a\overline{m}=0\rbrace $. Since $soc(A)$ is a homogeneous ideal, we define the {\it socle-vector} of $A$ as $s(A)=s=(s_0,s_1...,s_e)$, where $s_i=\dim_k soc(A)_i$. Note that $h_0=1$, $s_0=0$ and $s_e=h_e>0$.\\\indent
An algebra $A$ having a socle-vector $s$ of the form $s=(0,0,...,0,s_e)$ is called {\it level}. In particular, if $s_e=1$, $A$ is {\it Gorenstein}.\\\indent
Let us now recall the main facts of the theory of Inverse Systems which we will use throughout the paper. For a complete introduction, we refer the reader to $[Ge]$ and $[IK]$.\\\indent 
Let $S=k[y_1,...,y_r]$, and consider $S$ as a graded $R$-module where the action of $x_i$ on $S$ is partial differentiation with respect to $y_i$.\\\indent 
There is a one-to-one correspondence between artinian algebras $R/I$ and finitely generated $R$-submodules $M$ of $S$, where $I=M^{-1}$ is the annihilator of $M$ in $R$ and, conversely, $M$ is the $R$-submodule of $S$ which is annihilated by $I$ (cf. $[Ge]$, Remark 1), p. 17).\\\indent 
If $R/I$ has socle-vector $s$, then $M$ is minimally generated by $s_i$ elements of degree $i$, for $i=1,...,e$, and the $h$-vector of $R/I$ is given by the number of linearly independent derivatives in each degree obtained by differentiating the generators of $M$ (cf. $[Ge]$, Remark 2), p. 17).\\\indent
Remember that the {\it (degree) lexicographic order} is a total ordering, briefly indicated by $\lq \lq >"$, on the monomials of $R$ such that, if $M=x_1^{m_1}\cdot \cdot \cdot x_r^{m_r}$ and $N=x_1^{n_1}\cdot \cdot \cdot x_r^{n_r}$ have the same degree, then $M>N$ if and only if the first non-zero difference $m_i-n_i$ is positive.\\\indent
A monomial ideal $I$ of $R$ is a {\it lex-segment ideal} if, for every degree $d$ and monomials $T,T^{'}\in R_d$, whenever $T\in I_d$ and $T^{'}>T$, then $T^{'}\in I_d$.\\\indent
It is easy to see that the Inverse System module $M$ of a monomial ideal $I$ is generated in each degree (as a $k$-vector space) by monomials, precisely those of $S$ except for the monomials that (written in the $x_i$'s) belong to $I$. In particular, if $I$ is a lex-segment ideal, then $M$ is generated in each degree $d$ (as a $k$-vector space) by the last $h_d=\dim_k (R/I)_d$ monomials in the lexicographic order.\\
\\\indent
{\bf Definition-Remark 1.1.} Let $n$ and $i$ be positive integers. The {\it i-binomial expansion of n} is $$n_{(i)}={n_i\choose i}+{n_{i-1}\choose i-1}+...+{n_j\choose j},$$ where $n_i>n_{i-1}>...>n_j\geq j\geq 1$.\\\indent
Under these hypotheses, the $i$-binomial expansion of $n$ is unique (e.g., see $[BH]$, Lemma 4.2.6).\\\indent
Following $[BG]$, define, for any integer $a$, $$(n_{(i)})_{a}^{a}={n_i+a\choose i+a}+{n_{i-1}+a\choose i-1+a}+...+{n_j+a\choose j+a}.$$
\\\indent
{\bf Theorem 1.2} ($[BH]$, Theorem 4.2.10). {\it Let $h=(h_i)_{i\geq 0}$ be a sequence of non-negative integers, such that $h_0=1$, $h_1=r$ and $h_i=0$ for $i>e$. Then $h$ is the $h$-vector of some standard graded artinian algebra if and only if, for every $d$, $1\leq d\leq e-1$, $$h_{d+1}\leq ((h_d)_{(d)})_{+1}^{+1}.$$}
\\\indent
We just remark that Theorem 1.2 is due to Macaulay, and that it holds, with appropriate modifications, for any standard graded algebra, not necessarily artinian.\\
\\\indent
{\bf Lemma 1.3} ($[BG]$, Lemma 3.3). {\it Let $a,b$ be positive integers, $b>1$. Then the smallest integer $c$ such that $a\leq (c_{(b-1)})_{+1}^{+1}$ is $$c=(a_{(b)})_{-1}^{-1}.$$}
\\\indent
{\bf Remark 1.4.} This result supplies a lower-bound for the $i$-th entry of an $h$-vector, once the $(i+1)$-st entry is known. In terms of Inverse Systems, it supplies a lower-bound for the number of linearly independent first derivatives of any given set of linearly independent forms of degree $i+1$.\\
\\\indent
A sequence of non-negative integers which satisfies the growth condition of Macaulay's theorem is called an {\it $O$-sequence}.\\\indent
We place a partial ordering on tuples of integers of the same length, by defining $t\geq t^{'}$ if every entry of $t$ is greater than or equal to the corresponding entry of $t^{'}$.\\\indent
Given a socle-vector $s$, we say that an $h$-vector $h$ is the {\it minimum $h$-vector for $s$} (if it exists), if it is the (entry by entry) minimum among the $h$-vectors of all artinian algebras having $s$ as a socle-vector. In particular, the entry of degree 1 of such an $h$ is called the {\it minimum codimension for $s$} and is the least codimension that an algebra with socle-vector $s$ may have.\\\indent
Similarly, given an $h$-vector $h$, we say that a socle-vector $s$ is the {\it maximum socle-vector for $h$} (if it exists), if it is the (entry by entry) maximum among the socle-vectors of all artinian algebras having $h$ as an $h$-vector.\\
\\\indent
In the next section we construct a bijection between the set of all the finite $O$-sequences and the set of all the non-zero finite tuples of non-negative integers whose first entry is equal to 0. It is not difficult to show that the latter set coincides with the set of all the socle-vectors of standard graded artinian algebras. This bijection associates to a socle-vector $s$ its minimum $h$-vector $h$, which is shown to exist and is found explicitly, and to this $h$ its maximum socle-vector, which is shown to exist and to be equal to the $s$ above. In particular, for any socle-vector $s$, we find its minimum codimension. Our construction extends some results of $[BG]$, where Bigatti and Geramita studied level algebras.\\\indent
In Sections 3 to 5, we consider the main problem of this article: determining when there is a unique socle-vector associated to a given $h$-vector. In Section 3, we use some homological techniques to determine when, for a given $h$-vector $h$ of codimension $r$, the last module of any MFR of an algebra with $h$-vector $h$ is fixed. For these $h$'s, the associated socle-vector must then be unique.\\\indent
Then, in Section 5, we make a suitable use of Inverse Systems and lex-segment ideals to show that, (in most cases) where the maximum socle-vector for $h$ has two non-zero consecutive entries, there must be more than one socle-vector associated to $h$.\\\indent
In the case of codimension $r\leq 3$, we are able to solve the problem completely, by characterizing the $h$-vectors that admit a unique socle-vector (Section 4).\\\indent
The results obtained in this paper are part of the author's Ph.D. dissertation, written at Queen's University (Kingston, Ontario, Canada), under the supervision of Professor A.V. Geramita.\\
\\
\section{A bijection}
\indent

{\bf Theorem 2.1.} {\it i). Let $s=(0,s_1,...,s_e)$ be a socle-vector. Then the minimum $h$-vector for $s$ is $h=(1,h_1,...,h_e),$ where $h_e=s_e$, and inductively, for $i=1,...,e-1$, $$h_i=((h_{i+1})_{(i+1)})_{-1}^{-1}+s_i.$$
\indent In particular the minimum codimension for $s$ is $$r=h_1=((...((((s_e)_{(e)})_{-1}^{-1}+s_{i-1})_{(e-1)})_{-1}^{-1}+s_{e-2}...)_{(2)})_{-1}^{-1}+s_1.$$\indent
ii). Let $h=(1,h_1,...,h_e)$ be an $h$-vector. Then the maximum socle-vector for $h$ is $s=(0,s_1,...,s_e)$, where $s_e=h_e$, and inductively, for $i=1,...,e-1$, $$s_i=h_i-((h_{i+1})_{(i+1)})_{-1}^{-1}.$$\indent
iii). $h$ is the minimum $h$-vector for the socle-vector $s$ if and only if $s$ is the maximum socle-vector for the $h$-vector $h$.\\}
\\\indent
{\bf Proof.} i). Of course $h_e=s_e$. By Inverse Systems and Lemma 1.3, arguing by induction, given $h_{i+1}$ linearly independent forms of degree $i+1$, the least possible number of derivatives supplied in degree $i$ is $((h_{i+1})_{(i+1)})_{-1}^{-1}$. When we add $s_i$ linearly independent forms in degree $i$, we obtain the desired value for $h_i$. This shows that the $h$-vector of the statement is a lower-bound.\\\indent
In order to actually achieve it, we construct the following Inverse Systems module: consider, according to the lexicographic order, the last $h_e=s_e$ monomials of degree $e$ in $S=k[y_1,...,y_r]$, where $r=h_1$. Differentiating them, we obtain, without considering the coefficients, exactly the last $((h_{e})_{(e)})_{-1}^{-1}$ monomials of degree $e-1$, and then, adding to these the immediately previous $s_{e-1}$ monomials, we get the last $h_{e-1}$ monomials of degree $e-1$. By induction, we repeat the same process until we obtain all the $h_1$ monomials of degree 1. Clearly, the $h$-vector given by the ideal (which is a lex-segment ideal) annihilating the Inverse System module constructed in this way is the lower-bound of the statement.\\\indent
ii). Of course $s_e=h_e$. As we observed above, by Inverse Systems, the least possible number of derivatives given in degree $i$ by $h_{i+1}$ linearly independent forms of degree $i+1$ is $((h_{i+1})_{(i+1)})_{-1}^{-1}$. Hence, the Inverse System module has at most $h_i-((h_{i+1})_{(i+1)})_{-1}^{-1}$ new generators in degree $i$. Using induction, this suffices to prove the upper-bound for $s$.\\\indent
By Inverse Systems, reasoning as above, it is easy to see that, if $I$ is a lex-segment ideal and $R/I$ has $h$-vector $h$, then $R/I$ has the upper-bound $s$ of the statement as its socle-vector. This completes the proof of ii).\\\indent
iii). It immediately follows from i) and ii).{\ }{\ }\qed \\
\\\indent
{\bf Remark 2.2.} i). The result of Theorem 2.1, ii) was already known to Bigatti, Hulett and Pardue, who proved, more generally, that the lex-segment ideal supplies the maximal MFR associated to a given $h$-vector: see $[Bi]$, $[Hu]$ and $[Pa]$. We will state this theorem in the next section.\\\indent
ii). Theorem 2.1, iii) generalizes Proposition 3.7 and Theorem 5.8 of $[BG]$ to arbitrary socle-vectors.\\\indent
iii). It follows from Theorem 2.1 that there is a bijection between the set of the socle-vectors and the set of the $h$-vectors of standard graded artinian algebras, given by associating to a socle-vector its minimum $h$-vector.\\
\\\indent
{\bf Example 2.3.} Let $s=(0,0,1,0,2,4,2,5).$ Then $e=7$; $h_7=s_7=5$; $h_6=(5_{(7)})_{-1}^{-1}+2=5+2=7$; $h_5=(7_{(6)})_{-1}^{-1}+4=6+4=10$; $h_4=(10_{(5)})_{-1}^{-1}+2=9+2=11$; $h_3=(11_{(4)})_{-1}^{-1}+0=9+0=9$; $h_2=(9_{(3)})_{-1}^{-1}+1=6+1=7$; $h_1=(7_{(2)})_{-1}^{-1}+0=4+0=4$.\\\indent
Therefore, the minimum $h$-vector for the socle-vector $s$ is $$h=(1,4,7,9,11,10,7,5).$$\indent
In particular, the minimum codimension for $s$ is $r=h_1=4$.\\\indent
A similar computation shows that, as we expect from Theorem 2.1, iii), the maximum socle-vector for $h=(1,4,7,9,11,10,7,5)$ is $s=(0,0,1,0,2,4,2,5).$\\
\\
\section{Uniqueness of the socle-vector}
\indent

In the previous section we saw how to determine the maximum socle-vector for a given $h$-vector. In this section we study the problem of determining the $h$-vectors $h$, of arbitrary codimension $r$, that admit a unique socle-vector, i.e. such that all the algebras having $h$-vector $h$ have the maximum socle-vector for $h$ as their socle-vector. The simplest case is that of the minimum $h$-vector $h$ for a level socle-vector $s$: in fact, by $[BG]$, Proposition 3.7 (or by our Theorem 2.1, iii) above), we have that $s$ must be the only socle-vector associated to $h$.\\
\\\indent
We first recall some basic definitions and results about resolutions. As above, let $R=k[x_1,...,x_r]$ and let $I$ be a homogeneous ideal of $R$ such that $A=R/I$ is an artinian algebra, having $h$-vector $h=(h_0=1,h_1=r,h_2,...,h_e)$. The {\it minimal free resolution} ({\it MFR}, in brief) of $A$ is an exact sequence of $R$-modules of the form: $$0\longrightarrow F_r\longrightarrow F_{r-1}\longrightarrow ...\longrightarrow F_1\longrightarrow R\longrightarrow R/I \longrightarrow 0,$$
where, for $i=1,...,r$, $$F_i=\bigoplus_{j=1}^{n_i}R^{\beta_{i,j}}(-j),$$  and all the homomorphisms have degree 0.\\\indent The $\beta_{i,j}$'s are called the {\it graded Betti numbers} of $A$.\\\indent 
Then $\beta_{1,j}$ is the number of generators of $I$ in degree $j$. It is well-known that $F_r=\oplus_{j=1}^eR^{s_j}(-j-r)\neq 0,$ where $s=(0,s_1,...,s_e)$ is the socle-vector of $A$. Hence, the socle-vector may also be computed by considering the graded Betti numbers of the last module of the MFR of $A$.\\\indent  
With $h$ as above, let $h(z)=\sum_{i=0}^eh_iz^i$. MFR and $h$-vector of $A$ are related by the following well-known formula (see, e.g., $[FL]$, p. 131, point (j) for a proof): \begin{equation}\label{sh}h(z)(1-z)^r=1+\sum_{i,j}(-1)^i\beta_{i,j}z^j.\end{equation}
\\\indent
{\bf Theorem 3.1} ($[Bi]$,$[Hu]$,$[Pa]$). {\it Let $I$ and $J$ be ideals of $R$ such that $I$ is a lex-segment ideal and $R/I$ and $R/J$ are artinian algebras having the same $h$-vector. Let $\beta_{i,j}^I$ and $\beta_{i,j}^J$ be the graded Betti numbers of $I$ and $J$ respectively. Then $\beta_{i,j}^J\leq \beta_{i,j}^I$ for every $i$ and $j$.}\\
\\\indent
Following $[EK]$, let $T$ be a monomial in $R$, and define $m(T)=\max \lbrace i {\ } \mid {\ } x_i$ divides $T\rbrace $.\\
\\\indent
{\bf Theorem 3.2} ($[EK]$). {\it Let $I$ be a lex-segment ideal of $R$ such that $R/I$ has graded Betti numbers $\beta_{i,j}$. Denote by $G(I)$ the set of minimal (monomial) generators of $I$ and by $G(I)_d$ the subset of those of degree $d$. Then $$\beta_{i,j}=\sum_{T\in G(I)_{j-i+1}}{m(T)-1\choose i-1}.$$}
\\\indent
It is easy to see that, by equation (\ref{sh}) and Theorem 3.1, the MFR of any algebra having a given $h$-vector $h$ is obtained by cancelations from the MFR given by the lex-segment ideal associated to $h$. That is, for each shift $d$, the alternating sum of the exponents of the summands $R(-d)$ appearing in each module of the MFR must be constant for all the algebras having the same $h$-vector. The following result of Peeva shows that, moreover, these cancelations cannot be arbitrary:\\  
\\\indent
{\bf Theorem 3.3} ($[Pe]$). {\it Every cancelation must occur in adjacent modules of the MFR.}\\
\\\indent
{\bf Definition-Remark 3.4.} We say that an $h$-vector $h$ {\it admits a possible cancelation} if, in the MFR given by the lex-segment ideal associated to $h$, there is a common summand $R(-d)$ in the last and next to last module. Otherwise, we say that $h$ {\it admits no possible cancelations}.\\\indent
Notice that in this definition we make no comments about cancelations in other points of the MFR. We use the phrase \lq \lq no possible cancelations" just because our interest is only in the socle and hence in the last module of the MFR.\\
\\\indent
The following theorem is a very useful improvement of a result of $[CI]$ and $[GHMS]$.\\
\\\indent
{\bf Theorem 3.5.} {\it Let $h_{d-1}^{'}$, $h_{d}^{'}$ and $h_{d+1}^{'}$ be three integers such that $((h_{d}^{'})_{(d)})_{-1}^{-1}=h_{d-1}^{'}$ and $((h_{d}^{'})_{(d)})_{+1}^{+1}=h_{d+1}^{'}$. Suppose that $h_{d-1}^{'}+\alpha $, $h_{d}^{'}$ and $h_{d+1}^{'}$, for some integer $\alpha \geq 0$, are the entries of degree $d-1$, $d$ and $d+1$ of the $h$-vector of an algebra $A$. Then the entry of degree $d-1$ of the socle-vector of $A$ is equal to $\alpha $.}\\
\\\indent
{\bf Proof.} Let $h$ be the $h$-vector of $A$, $I$ the lex-segment ideal associated to $h$, and $\beta_{i,j}$ the Betti numbers of $R/I$. By Theorems 3.1 and 3.3 and the relationship between the socle-vector of an algebra and the last module of its MFR, it is clearly enough to show that $\beta_{r,d-1+r}=\alpha $ and $\beta_{r-1,d-1+r}=0$.\\\indent
We saw in the proof of Theorem 2.1 that the condition $((h_{d}^{'})_{(d)})_{-1}^{-1}=h_{d-1}^{'}$ implies that the entry of degree $d-1$ of the socle-vector of $R/I$ is equal to $\alpha $, i.e. $\beta_{r,d-1+r}=\alpha $. Moreover, from the maximal growth $((h_{d}^{'})_{(d)})_{+1}^{+1}=h_{d+1}^{'}$, we have that $I$ has no generators in degree $d+1$. Hence, by Theorem 3.2, $\beta_{r-1,d-1+r}=0$. This proves the theorem.{\ }{\ }\qed \\
\\\indent
{\bf Example 3.6.} Our Theorem 3.5 improves $[GHMS]$, Theorem 3.4 (which was a result first stated in $[CI]$, Remark 2.8). This is true since $((h_{d-1}^{'})_{(d-1)})_{+1}^{+1}=h_{d}^{'}$ implies $((h_{d}^{'})_{(d)})_{-1}^{-1}=h_{d-1}^{'}$, but not vice versa (just an inequality holds).\\\indent
For example, let $h=(1,4,8,9,12,13)$. Since $((9)_{(3)})_{-1}^{-1}=6$ (even if $((6)_{(2)})_{1}^{1}=10>9$) and $((9)_{(3)})_{1}^{1}=12$, Theorem 3.5 implies that any algebra having $h$-vector $h$ must have the entry of degree 2 of its socle-vector equal to 2; instead, $[GHMS]$, Theorem 3.4 cannot be applied to this example.\\
\\\indent
{\bf Corollary 3.7.} {\it Let $h=(1,h_1=r,h_2,...,h_e)$ be an $h$-vector such that, for every $i\geq 2$, the inequality $((h_{i})_{(i)})_{-1}^{-1}\neq h_{i-1}$ implies the equality $((h_{i})_{(i)})_{+1}^{+1}=h_{i+1}$. Then there exists a unique socle-vector associated to $h$.}\\
\\\indent
{\bf Proof.} It suffices to apply the same argument of Theorem 3.5 to any non-zero entry of the socle-vector of $R/I$, where $I$ is the lex-segment ideal that gives the $h$-vector $h$, in order to obtain that $h$ admits no possible cancelations. In particular, the socle-vector associated to $h$ must be unique.{\ }{\ }\qed \\
\\\indent
{\bf Remark 3.8.} In Corollary 3.7, we have shown that the hypothesis on $h$ implies that $h$ admits no possible cancelations. If $r\leq 3$, the converse is also true. In other words, for $r\leq 3$, the hypothesis on $h$ of Corollary 3.7 is equivalent to the fact that $h$ admits no possible cancelations.\\\indent
Indeed, by Theorem 3.2 it is easy to see that, for any $r$, $h$ admits no possible cancelations if and only if, for every $i>2$ such that $((h_{i})_{(i)})_{-1}^{-1}\neq h_{i-1}$ (i.e. $s_{i-1}\neq 0$), both $y_r$ and $y_{r-1}$ do not divide any generator of $I$ of degree $i+1$ (i.e. $\beta_{r-1,i-1+r}=0$).\\\indent
In particular, for $r\leq 3$, $\beta_{r-1,i-1+r}=0$ if and only if $I$ has no generators in degree $i+1$, if and only if $((h_{i})_{(i)})_{+1}^{+1}=h_{i+1}$. The result now follows.\\
\\
\section{Characterization in codimension $r\leq 3$}
\indent

As we saw in Remark 3.8, for $r\leq 3$ it is possible to characterize the $h$-vectors $h$ that admit a possible cancelation (see Definition-Remark 3.4). Now we show that the \lq \lq best-case\rq \rq {\ }scenario occurs, i.e. any time $h$ admits a possible cancelation, a cancelation actually happens and, therefore, there is more than one socle-vector associated to $h$.\\\indent
Notice that the fact that one cancelation actually occurs (i.e. that one summand $R(-d)$ cancels) does not imply the much stronger fact that {\it all} the common summands cancel from the last and next to last module of the MFR given by the lex-segment ideal associated to $h$. (See $[GHMS]$ for examples of $h$-vectors for which no level algebra exists, even though the MFR of the lex-segment ideal would allow the possibility of performing all the necessary cancelations.)\\
\\\indent
{\bf Theorem 4.1.} {\it Let $h=(1,h_1,...,h_e)$ be an $h$-vector, with $h_1=r\leq 3$. Then $h$ admits a possible cancelation if and only if a cancelation actually occurs.}\\
\\\indent
From Remark 3.8 and Theorem 4.1, we have at once:\\
\\\indent 
{\bf Corollary 4.2.} {\it Let $h=(1,h_1=r,h_2,...,h_e)$, with $r\leq 3$. Then there exists a unique socle-vector associated to $h$ if and only if, for every $i\geq 2$, the inequality $((h_{i})_{(i)})_{-1}^{-1}\neq h_{i-1}$ implies the equality $((h_{i})_{(i)})_{+1}^{+1}=h_{i+1}$.}\\
\\\indent
{\bf Proof of Theorem 4.1.} Of course, the only implication that needs to be shown is \lq \lq $\Longrightarrow $\rq \rq . Let $r=3$, and let $I\subseteq R=k[x,y,z]$ be the lex-segment ideal associated to $h$. Suppose $h$ admits a possible cancelation. Then, by Remark 3.8, there is an integer $d$, $4\leq d<e+3$, such that we have the inequalities $((h_{d-2})_{(d-2)})_{-1}^{-1}<h_{d-3}$ and $((h_{d-2})_{(d-2)})_{+1}^{+1}>h_{d-1}$, which together say the possible cancelation in the last module of the MFR of $R/I$ is of an $R(-d)$. In particular, $R/I$ has a non-zero socle in degree $d-3$. Thus, in order to prove that a cancelation actually occurs, it is enough to find an ideal $J$ of $R$ such that $R/J$ has the same $h$-vector $h$ as $R/I$ but a smaller socle in degree $d-3$. We now seek to construct such a $J$.\\\indent
Let $I_{d-2}$ be the lex-segment $<T_1=x^{d-2},T_2,...,T_p>$. At least one monomial in $I_{d-2}$ is divisible by $z$ (since  $R/I$ has a non-zero socle in degree $d-3$), and let $T_q$ be the last such monomial. Also, it is easy to see that, in a lex-segment of monomials of a given degree, there are never two consecutive terms not divisible by $z$ (except for the first two). Hence $q=p$ or $q=p-1$. First suppose $q=p$.\\\indent
Under this hypothesis we define $J$ as follows. Let $J_i=I_i$ for $i<d-2$.\\\indent
We now construct $J_{d-2}$. Since $((h_{d-2})_{(d-2)})_{-1}^{-1}<h_{d-3}$, we have $((h_{d-3})_{(d-3)})_{+1}^{+1}>h_{d-2}$, and thus $I$ has generators in degree $d-2$, one of them being necessarily $T_p$. Let $J_{d-2}=<T_1,...,T_{p-1},\overbrace{T_p},T_{p+1}>$, where the over braced term $T_p$ is omitted, i.e. $J_{d-2}$ is spanned by the same monomials as $I_{d-2}$, except for $T_p$ which is replaced by $T_{p+1}$.\\\indent
Notice that $T_p/z\notin I_{d-3}=J_{d-3}$ (since $T_p$ is a generator of $I$), and that the class of $T_p/z$ (which is then non-zero) is, by the properties of the lexicographic order, in the socle of $R/I$, while, by construction, the class of $T_p/z$ is not in the socle of $R/J$.\\\indent
Thus, in order to show that the socle of $R/J$ is smaller than that of $R/I$ in degree $d-3$, it suffices to prove that all the monomials that give elements of the socle of $R/J$ in degree $d-3$ give also elements of the socle of $R/I$. If a monomial $M\in R_{d-3}$ gives an element of the socle of $R/J$ and not an element of the socle of $R/I$, then clearly $Mz=T_{p+1}$. I.e., it is enough to show that either $T_{p+1}$ is not divisible by $z$ or, if it is, then the class of $T_{p+1}/z$ is not in the socle of $R/J$.\\\indent
Write $T_p=x^ay^bz^c$, with $c>0$. If $b=0$, then $T_{p+1}=x^{a-1}y^{c+1}$ (which is not divisible by $z$). If $b>0$, then $T_{p+1}=x^ay^{b-1}z^{c+1}$. Thus, $(T_{p+1}/z)y=T_p$, and therefore  the class of $T_{p+1}/z$ does not belong to the socle of $R/J$, as we wanted to show.\\\indent
At this point we have partially constructed an ideal $J$ (up to degree $d-2$) so that $R/J$ has the same $h$-vector of $R/I$ and (up to degree $d-3$) $R/J$ has socle strictly smaller that that of $R/I$. So, to prove the theorem (for $q=p$), it is enough  to show that we can extend $J$ to an ideal of $R$ such that the $h$-vector of $R/J$ is that of $R/I$.\\\indent
We now want to construct $J_{d-1}$. From the inequality $((h_{d-2})_{(d-2)})_{+1}^{+1}>h_{d-1}$, we see that $I$ must have generators in degree $d-1$. Note that, by construction, $R_1J_{d-2}=<T_1,...,T_{p-1}z,\overbrace{?,T_pz},?,T_{p+1}z>$, where \lq \lq ?\rq \rq {\ }denotes terms (if they exist) between two monomials divisible by $z$.\\\indent
{\it Case 1.} $I$ has more than one generator in degree $d-1$. In this case, one of these generators is necessarily
$T_{p+1}z$ (because $T_pz\in R_1I_{d-2}$). Since $T_{p+1}z$ is the last monomial in $R_1J_{d-2}$, then $R_1J_{d-2}\subseteq I_{d-1}$. So, it suffices to add to $R_1J_{d-2}$ enough monomials to make it equal to $I_{d-1}$. Thus, in this case, we can take $J_i=I_i$ for all $i\geq d-1$, and the construction is complete.\\\indent
{\it Case 2.} $I$ has only one generator in degree $d-1$. If the last \lq \lq?\rq \rq {\ }in the terms spanning $R_1J_{d-2}$ does not exist, then, by the same reasoning, we are done. Therefore suppose it exists, and call it $T$. Thus, $T$ is the only generator of $I$ having degree $d-1$, and since, in the lexicographic order, it is between $T_pz$ and $T_{p+1}z$, we must have $T=T_{p+1}y$. By definition, $z$ does not divide $T$, and hence does not divide $T_{p+1}$.\\\indent
Let $J_{d-1}=<T_1,...,T_py,\overbrace{T_pz},T,T_{p+1}z>$ (which, of course, contains $R_1J_{d-2}$). Hence $J_{d-1}$ has the same dimension as $I_{d-1}$, since they are equal except for $T_pz$ which is replaced by $T_{p+1}z$ in $J_{d-1}$.\\\indent
For all $i\geq d$, we want to construct $J_i$ having the same dimension as $I_i$. If $I$ has generators in degree $d$, then we are done (by the same argument as above, since there is only one monomial in $R_1J_{d-1}$ not in $R_1I_{d-1}$, namely $T_{p+1}z^2$, so we only need to add it in order to get $J_d=I_d$, and that is enough). If $I$ has no generators in degree $d$, then add no generators to $J$ in degree $d$ either, since the dimensions of $J_d=R_1J_{d-1}$ and $I_d=R_1I_{d-1}$ are equal. In fact, it is easy to see that $J_d=I_d$, except that $T_{p+1}z^2$ belongs to $J_d$ but not to $I_d$, and $T_{p}z^2$ is in $I_d$ but not in $J_d$ (since $z$ divides $T_p$ but not $T_{p+1}$).\\\indent
By induction, we continue this procedure until, in some degree $i\geq d$, $I$ has generators, and then we have finished. Otherwise, $R_1J_{i-1}$ and $R_1I_{i-1}$ have the same dimension: in fact, similarly to above, they are equal, except that $T_{p+1}z^{i-d+2}$ belongs to $J_i$ but not to $I_i$, and $T_{p}z^{i-d+2}$ is in $I_i$ but not in $J_i$. Since $I$ necessarily has generators in degree $e+1$, this process must eventually end, and the construction is complete. This proves the case $q=p$.\\\indent
Now let $q=p-1$, i.e. $T_p$ is not divisible by $z$. Then $T_p$ has the form $T_p=x^{a-1}y^{c+1}$, and therefore $T_q=x^az^c$, with $c>0$. It follows that $T_{p+1}=x^{a-1}y^cz$ and $T_{p+2}=x^{a-1}y^{c-1}z^2$. As above, we want to construct an ideal $J$ of $R$ such that $R/J$ has the same $h$-vector $h$ as $R/I$ but a smaller socle in degree $d-3$. (Notice that now $R/J$ cannot be defined as we did before, since the class of $T_{p+1}/z$ would be in its socle and therefore the argument above would not apply.)\\\indent
We seek to construct $J$. First let $J_i=I_i$ for $i<d-2$. Notice that $T_q/z$ does not belong to $I_{d-3}$. In fact, any monomial that is after $T_q/z$ in the lexicographic order, multiplied by $z$ does not belong to $I_{d-2}$, and therefore, if $T_q/z\in I_{d-3}$, the socle of $R/I$ in degree $d-3$ would be zero, a contradiction. In particular, it follows that $T_q$ is a generator of $I$: in fact, $T_q/x$ and $T_q/y$, if they exist, are after $T_q/z$ in the lexicographic order, and therefore do not belong to $I_{d-3}$ either.\\\indent
{\it Case $c>1$.} Let $J_{d-2}=<T_1,...,T_{q-1},\overbrace{T_q},T_p,\overbrace{T_{p+1}},T_{p+2}>$. Reasoning similarly to the case $q=p$, we have that $T_q/z$ gives an element that belongs to the socle of $R/I$ in degree $d-3$ (where it is non-zero by the consideration above), but not to that of $R/J$, which is then smaller (since, moreover, $(T_{p+2}/z)y=T_{p+1}$, and therefore  the class of $T_{p+2}/z$ is not in the socle of $R/J$).\\\indent
Now $R_1J_{d-2}=<T_1x,...,T_qy,\overbrace{T_qz},T_py,...,T_{p+2}z>$ (which is a lex-segment minus $T_qz$). As above, if $I$ has more than one generator in degree $d-1$, we are done (because, since $T_p=x^{a-1}y^{c+1}$, we have $R_1J_{d-2}\subseteq I_{d-1}$). Otherwise, if $I$ has only one generator in degree $d-1$, it is easy to see that the same idea we used in the case $q=p$ also applies to construct the desired ideal $J$ in this case.\\\indent
{\it Case $c=1$.} Since $T_q=x^{d-3}z$ is a generator of $I$, we clearly have that $I_{d-3}=0$, whence $I_{d-2}=<x^{d-2},x^{d-3}y,x^{d-3}z,x^{d-4}y^2>$ and the socle of $R/I$ in degree $d-3$ is spanned by the class of $x^{d-3}$.\\\indent
Let $J_{d-2}=<x^{d-2},x^{d-3}z,x^{d-4}y^2,x^{d-4}yz>$. We have $soc(R/J)_{d-3}=0$ and $R_1J_{d-2}=<x^{d-1},...,x^{d-4}y^2z,x^{d-4}yz^2>$ (which is a lex-segment). Therefore $R_1J_{d-2}$ differs from $R_1I_{d-2}$ only in containing the monomial $x^{d-4}yz^2$. Thus, since $I$ has generators in degree $d-1$, we have $R_1J_{d-2}\subseteq I_{d-1}$. At this point, as we observed above, it suffices to add the right generators to $J$ in order to obtain $J_i=I_i$ for all $i\geq d-1$ and to complete the construction. This proves the theorem for $r=3$.\\\indent
Let $r=2$. The proof of the case $q=p$ for $r=3$ applies {\it mutatis mutandis} (except that the situation now is simpler). For this reason we will omit the complete argument and, with the notation above, just make some remarks about the main differences in the proof (none of which are substantial).\\\indent
If a cancelation in the last (i.e. in the second) module of the MFR of $k[x,y]/I$ is possible in shift $(-d)$, then $I$ has generators in degree $d$, and the socle of $R/I$ is non-zero in degree $d-2$. This time, let $J_i=I_i$ for all $i<d-1$. Let $T_p$ be the last monomial in $I_{d-1}$, and define $T_q$ as the last monomial of $I_{d-1}$ divisible by $y$; we always have $p=q$, since $y$ divides all the monomials after $x^{d-1}$. Define $J_{d-1}$ in a way analogous to the way we defined $J_{d-2}$ in the case $p=q$ above, and observe that $T_p/y$ does not belong to $I$, that $T_p/y$ gives an element which is non-zero in the socle of $R/I$ but which is not in the socle of $R/J$, and that nothing new is in the socle of $R/J$. The construction is quicker, since there is only one monomial in $R_1J_{d-1}$ but not in $R_1I_{d-1}$ (namely $T_{p+1}y$), but $I$ has generators in degree $d$, and therefore we can make $J_i=I_i$ for all $i\geq d$ as above.\\\indent
This completes the proof of the case $r=2$ and that of the theorem.{\ }{\ }\qed \\
\\\indent
{\bf Example 4.3.} Let $h=(1,3,6,10,12,14)$. By Theorem 2.1, the maximum socle-vector for $h$ is $s=(0,0,0,1,0,14)$. Let $I\subseteq R=k[x,y,z]$ be the lex-segment ideal associated to $h$. By Remark 3.8 (and Corollary 3.7), the MFR of $R/I$ admits a possible cancelation in the last module. More precisely, an $R(-6)$ can be canceled. We want to find, following the procedure suggested by the proof of Theorem 4.1, an algebra $R/J$ whose MFR is obtained from that of $R/I$ by performing that cancelation.\\\indent
We have $I_i=0$ for $i<4$, $I_4=<x^4,x^3y,x^3z>$, $I_5=<x^5,...,x^3z^2,x^2y^3>$ (with $x^2y^3$ a generator of $I$), $I_6=<x^6,...,x^2y^3z,x^2y^2z^2,...,z^6>=R_6$ (with the 16 monomials $x^2y^2z^2,...,z^6$ generators of $I$), and of course $I_i=R_i$ for $i>6$.\\\indent
Let us construct $J$ as follows: let $J_i=I_i=0$ for $i<4$, $J_4=<x^4,x^3y,x^2y^2>$, $J_5=R_1J_4=<x^5,...,x^3yz,x^2y^3,x^2y^2z>$, and $J_i=I_i=R_i$ for $i\geq 6$ (by adding the 16 monomials that are left to obtain all of $R_6$ from $R_1J_5$). Notice that $R/J$ has the same $h$-vector $h$ as $R/I$ but a zero-dimensional socle in degree 3, since our construction \lq \lq eliminated\rq \rq {\ }$x^3$, whose class was in the socle of $R/I$. Thus, $R/J$ is the (level) algebra we desired.\\
\\\indent
{\bf Remark 4.4.} i). In the proof of Theorem 4.1 we have actually shown more than the statement. In fact, we have proven that, every time $h$ admits a possible cancelation of some summand $R(-d)$ in the last module of the MFR given by its lex-segment ideal, then there exists an algebra for which this cancelation occurs, i.e. this algebra has a smaller socle in degree $d-3$ than that given by the lex-segment ideal.\\\indent
This implies that, if $h$ admits a possible cancelation in $n$ different shifts, then there exist (at least) $n+1$ different socle-vectors associated to $h$.\\\indent
ii). If $r=2$, the condition on $h$ so that $h$ admits no possible cancelations can be written more simply, i.e. Corollary 4.2 can be rephrased as follows: there exists a unique socle-vector associated to $h$ if and only if, for every $i\geq 2$, the inequality $h_i<h_{i-1}$ implies the equality $h_{i}=h_{i+1}$.\\
\\
\section{Non-uniqueness of the socle-vector}
\indent

In this section we continue the study of the $h$-vectors that admit more than one socle-vector. This time we work in arbitrary codimension $r$, where the situation is much more complicated. The main result is Theorem 5.4, which give some sufficient conditions for the existence of more than one socle-vector associated to a given $h$-vector $h$. We believe there are other cases where the uniqueness fails behind those considered here, and it would be of great interest to characterize the $h$-vectors admitting a unique socle-vector in any codimension.\\\indent
We first need some preliminary results.\\
\\\indent
{\bf Lemma 5.1} ($[Ia]$, Proposition 4.7). {\it Let $F=\sum_{t=1}^m L_t^d$ be a form of degree $d$ in $S=k[y_1,...,y_r]$, where the $L_t=\sum_{k=1}^rb_{tk}y_k$ are linear forms, and let $I_{\bf b}=Ann(F)\subseteq R$, with ${\bf b}=(b_{11},...,b_{mr})$. Then there exists a non-empty open subset $U$ of $k^{mr}$ such that, for every ${\bf b}\in U$, the Gorenstein artinian algebras $R/I_{\bf b}$ all have the same $h$-vector, denoted by: $$h(m,d)=(1,h_1(m,d),...,h_d(m,d)=1),$$ where, for $j=1,...,d$, $$h_j(m,d)=\min \lbrace m,\dim_kR_j,\dim_kR_{d-j}\rbrace .$$}
\\\indent 
{\bf Lemma 5.2} ($[Ia]$, Theorem 4.8 A). {\it Let $h=(1,h_1,...,h_d)$ be the $h$-vector of an algebra $A=R/I$, where $I$ annihilates the $R$-submodule $M$ of $S$. Let $m\leq {r-1+d\choose d}-h_d$. Then, if $F$ is the sum of the $d$-th powers of $m$ generic linear forms (the Gorenstein $h$-vector of $R/<F>^{-1}$ is $h(m,d)$, given by Lemma 5.1), then the algebra associated to $M^{'}=<M,F>$ has $h$-vector $H=(1,H_1,...,H_d),$ where, for $i=1,...,d$, $$H_i=\min \lbrace h_i+h_i(m,d),{r-1+i\choose i}\rbrace .$$}
\\\indent 
{\bf Lemma 5.3.} {\it Let $i$ and $a$ be positive integers, and let $(a_{(i)})_{-1}^{-1}=b$. Then:\\\indent
i). $(a_{(i)})_{+1}^{+1}<((a+1)_{(i)})_{+1}^{+1}$.\\\indent
ii). If $a\geq 2$, then $((a-1)_{(i)})_{-1}^{-1}\in \lbrace b,b-1\rbrace $.\\\indent
iii). If $a\geq 3$, then $((a-2)_{(i)})_{-1}^{-1}\in \lbrace b,b-1,b-2\rbrace $.}\\
\\\indent
{\bf Proof.} This is an easy exercise which is left to the reader.{\ }{\ }\qed \\
\\\indent
{\bf Theorem 5.4.} {\it Let $h=(1,h_1=r,h_2,...,h_e)$ be any $h$-vector, and let $s=(0,s_1,...,s_e)$ be the maximum socle-vector for $h$ (see Theorem 2.1). Define $t$ as the largest integer such that $h$ is generic (i.e. $h_t={r-1+t\choose t}$). If any one of the following conditions is verified, then there exists more than one socle-vector associated to $h$:\\\indent
i). $$s_{t+1}s_{t+2}\neq 0;$$\indent
ii). There exists an index $i$ such that $$s_is_{i+1}\neq 0 {\ }\text{ and }{\ } ((h_{i+1}+1)_{(i+1)})_{-1}^{-1}>((h_{i+1})_{(i+1)})_{-1}^{-1};$$\indent
iii). $$s_{e-1}\neq 0.$$}
\\\indent
{\bf Proof.} i). We know that the maximum socle-vector $s$ for $h$ is given by the algebra associated to the lex-segment ideal $I$, say with Inverse System $M$. We want to construct an Inverse System module $M^{'}$ which gives the $h$-vector $h$ and a socle-vector $s^{'}<s$.\\\indent
Let $M^{''}$ be the Inverse System module of another lex segment ideal having the same generators as $M$ in degrees higher than $t+2$, $s_{t+2}-1$ generators in degree $t+2$, $s_{t+1}-1$ in degree $t+1$ and $s_t$ in degree $t$ (they are uniquely determined, of course, as the last possible according to the lexicographic order; see the considerations in the Introduction).\\\indent
Since $((h_{t+2})_{(t+2)})_{-1}^{-1}=h_{t+1}-s_{t+1}$, by Lemma 5.3, ii) we have that $$((h_{t+2}-1)_{(t+2)})_{-1}^{-1}\in \lbrace h_{t+1}-s_{t+1}, h_{t+1}-s_{t+1}-1\rbrace .$$\indent
Consider a form $F$ of degree $t+2$. If $((h_{t+2}-1)_{(t+2)})_{-1}^{-1}=h_{t+1}-s_{t+1}$, let $F$ be the $e$-th power of one generic linear form, while, if $((h_{t+2}-1)_{(t+2)})_{-1}^{-1}=h_{t+1}-s_{t+1}-1$, let $F$ be the sum of the $e$-th powers of two generic linear forms. Define $M^{'}=<M^{''},F>$. By Lemmata 5.1, 5.2 and 5.3, it is easy to check that $M^{'}$ is the desired module.\\\indent
ii). Let $M$ be defined as above. Also in this case, we want to find a module $M^{'}$ which gives the $h$-vector $h$ and a socle-vector $s^{'}<s$.\\\indent
Let $M^{'}$ have the same generators of $M$ in all degrees lower than $i$ and higher than $i+1$, the first $s_i-1$ generators of $M$ according to the lexicographic order in degree $i$, and the same generators of $M$ in degree $i+1$ except for the first, which is replaced by the immediately previous monomial in the lexicographic order. We want to show that this module $M^{'}$ is the one we are looking for, i.e. that it gives the $h$-vector $h$.\\\indent
By Lemma 5.3, ii), the hypothesis $((h_{i+1}+1)_{(i+1)})_{-1}^{-1}>((h_{i+1})_{(i+1)})_{-1}^{-1}$ implies that $((h_{i+1}+1)_{(i+1)})_{-1}^{-1}=((h_{i+1})_{(i+1)})_{-1}^{-1}+1$, whence the last $h_{i+1}+1$ monomials of degree $i+1$ (call them, in decreasing order, $T$, $U$, $W$, ...) give, by differentiation (without considering the coefficients), the last $h_i-s_i+1$ monomials of degree $i$ (call them, always in decreasing order, $T^{'}$, $U^{'}$, $W^{'}$, ...).\\\indent
At this point, in order to show that $M^{'}$ gives the $h$-vector $h$, it clearly suffices to show that the monomials $T$, $W$, ... (the same as above after getting rid of $U$) give, by differentiation, all of $T^{'}$, $U^{'}$, $W^{'}$, ....\\\indent
Recall that, for a monomial $P\in k[y_1,...,y_r]$, $m(P)=\max \lbrace  i{\ } \mid {\ } y_i {\ }\text{divides}{\ } P\rbrace .$ Hence $T^{'}=T/y_{m(T)}$. Since $((h_{i+1}+1)_{(i+1)})_{-1}^{-1}=((h_{i+1})_{(i+1)})_{-1}^{-1}+1$, we have $U^{'}=U/y_{m(U)}$. Instead, $W/y_{m(W)}\in \lbrace W^{'},U^{'}\rbrace $. If $W/y_{m(W)}=U^{'}$, then we are done. Therefore, let us suppose, from now on, that $W/y_{m(W)}=W^{'}$. It is enough to show that $T/y_j=U^{'}$ for some $j<m(T)$.\\\indent
The equality $((h_{i+1}+1)_{(i+1)})_{-1}^{-1}=((h_{i+1})_{(i+1)})_{-1}^{-1}+1$ easily implies that $m(T)=r$ (this fact was first noted in $[CI]$). Similarly, by the supposition above, we have $m(U)=r$. Hence $T$ has the form $T=y_1^{\beta_1}\cdot \cdot \cdot y_l^{\beta_l}y_r^{\beta }=Qy_r^{\beta }$, for some $\beta_l,\beta >0$ and $l<r$. Therefore $U=(Q/y_l)y_{l+1}^{\beta +1 }$. Since $y_r$ divides $U$, we have $l+1=r$. Now it is easy to check that $T/y_{r-1}=U^{'}$, as we wanted to show.\\\indent
iii). Let $M$ be defined as above. Again, we want to find a module $M^{'}$ which gives the $h$-vector $h$ and a socle-vector $s^{'}<s$.\\\indent
Consider the last $h_{e}$ monomials of degree $e$ (call them, in decreasing order, $U$, $W$, ...). They give, by differentiation (always without considering the coefficients), the last $h_{e-1}-s_{e-1}$ monomials of degree $e-1$ (call them, always in decreasing order, $U^{'}$, $W^{'}$, ...). Let $T$ be the first monomial of degree $e$ before $U$, in the lexicographic order, which is divisible by $y_r$ (it must exist, otherwise one can easily see that $U^{'}=y_1^{e-1}$, and therefore $s_{e-1}=0$). Suppose that $((h_e+1)_{e})_{-1}^{-1}=((h_{e})_{(e)})_{-1}^{-1}$, otherwise we are in hypothesis ii). This means that between $T$ and $U$ there must be other monomials, say (from $T$ to $U$), $T_1$, $T_2$, ..., $T_p$. Furthermore, let $T^{'}$ be the monomial of degree $e-1$ immediately before $U^{'}$. Hence, $T/y_r=T^{'}$.\\\indent
It is enough to find $h_e$ monomials of degree $e$ that give, by differentiation, all of $T^{'}$, $U^{'}$, $W^{'}$, .... (Notice that now we have freedom in this choice, since there are no generators of $M$ in degrees higher than $e$.) In fact, if we take these $h_e$ monomials as generators of $M^{'}$ in degree $e$, in degree $e-1$ we choose the last $s_e-1$ generators of $M$, and the same generators as $M$ in all the lower degrees, then, by Inverse Systems, the $h$-vector given by $M^{'}$ is clearly the same as that given by $M$.\\\indent
Similarly to case ii), if $W/y_{m(W)}=U^{'}$, then we are done (it is enough to choose the last $h_e-1$ monomials of degree $e$ plus $T$). Therefore, let us suppose that $W/y_{m(W)}=W^{'}$. 
We have that $T$ has the form $T=y_1^{\beta_1}\cdot \cdot \cdot y_l^{\beta_l}y_r^{\beta }=Qy_r^{\beta }$, for some $\beta_l,\beta >0$ and $l<r$, whence $T_1=(Q/y_l)y_{l+1}^{\beta +1 }$. Moreover, $T^{'}=T/y_r=Qy_r^{\beta -1}$, and therefore $U^{'}=(Q/y_l)y_{l+1}^{\beta }$. This implies that $U^{'}=T_1/y_{l+1}$.\\\indent
At this point, it is easy to see that it suffices to choose, as generators of $M^{'}$ of degree $e$, $T$, $T_1$, and the last $h_e-1$ monomials in the lexicographic order except for $y_r^e$. This completes the proof of the theorem.{\ }{\ }\qed \\
\\\indent
{\bf Remark 5.5.} If $s$ satisfies condition i) of Theorem 5.4, the new socle-vector $s^{'}=(0,s_1^{'},...,s_e^{'})$ associated to $h$ that the proof of the theorem constructs is equal to $s$ except for $s_{t+1}^{'}=s_{t+1}-1$.\\\indent
If $s$ satisfies condition ii) of the theorem, then $s^{'}$ equals $s$ except for $s_i^{'}=s_i-1$.\\\indent
If $s$ satisfies condition iii) of the theorem, then $s^{'}$ equals $s$ except for $s_{e-1}^{'}=s_{e-1}-1$.\\
\\
\\
\\
{\bf \huge References}\\
\\
$[Bi]$ {\ } A. Bigatti: {\it Upper bounds for the Betti numbers of a given Hilbert function}, Comm. in Algebra 21 (1993), No. 7, 2317-2334.\\
$[BG]$ {\ } A.M. Bigatti and A.V. Geramita: {\it Level Algebras, Lex Segments and Minimal Hilbert Functions}, Comm. in Algebra 31 (2003), 1427-1451.\\
$[BH]$ {\ } W. Bruns and J. Herzog: {\it Cohen-Macaulay rings}, Cambridge studies in advanced mathematics, No. 39, Revised edition (1998), Cambridge, U.K..\\
$[CI]$ {\ } Y.H. Cho and A. Iarrobino: {\it Hilbert Functions and Level Algebras}, J. of Algebra 241 (2001), 745-758.\\
$[EK]$ {\ } S. Eliahou and M. Kervaire: {\it Minimal resolutions of some monomial ideals}, J. of Algebra 129 (1990), 1-25.\\
$[FL]$ {\ } R. Fr\"oberg and D. Laksov: {\it Compressed Algebras}, Conference on Complete Intersections in Acireale, Lecture Notes in Mathematics, No. 1092 (1984), 121-151, Springer-Verlag.\\
$[Ge]$ {\ } A.V. Geramita: {\it Inverse Systems of Fat Points: Waring's Problem, Secant Varieties and Veronese Varieties and Parametric Spaces of Gorenstein Ideals}, Queen's Papers in Pure and Applied Mathematics, No. 102, The Curves Seminar at Queen's (1996), Vol. X, 3-114.\\
$[GHMS]$ {\ } A.V. Geramita, T. Harima, J. Migliore and Y.S. Shin: {\it The Hilbert Function of a Level Algebra}, Memoirs of the Amer. Math. Soc., to appear.\\
$[Hu]$ {\ } H. Hulett: {\it Maximum Betti numbers of homogeneous ideals with a given Hilbert function}, Comm. in Algebra 21 (1993), No. 7, 2335-2350.\\
$[Ia]$ {\ } A. Iarrobino: {\it Compressed Algebras: Artin algebras having given socle degrees and maximal length}, Trans. Amer. Math. Soc. 285 (1984), 337-378.\\
$[IK]$ {\ } A. Iarrobino and V. Kanev: {\it Power Sums, Gorenstein Algebras, and Determinantal Loci}, Springer Lecture Notes in Mathematics (1999), No. 1721, Springer, Heidelberg.\\
$[Pa]$ {\ } K. Pardue: {\it Deformation classes of graded modules and maximal Betti numbers}, Illinois J. Math. 40 (1996), 564-585.\\
$[Pe]$ {\ } I. Peeva: {\it Consecutive cancellations in Betti numbers}, Proc. Amer. Math. Soc. 132 (2004), 3503-3507.

}

\end{document}